\documentclass[10pt]{article}

\usepackage[english]{babel}
\usepackage[utf8x]{inputenc}
\usepackage[T1]{fontenc}

\usepackage[letterpaper,top=1in,bottom=1in,left=1in,right=1in,marginparwidth=1.75cm]{geometry}

\usepackage{amsmath,amsfonts,amsthm,amssymb,mathrsfs,dsfont} 
\usepackage{mathtools}
\usepackage{graphicx}
\usepackage[colorinlistoftodos]{todonotes}
\usepackage[colorlinks=true, allcolors=blue]{hyperref}
\usepackage[capitalize,nameinlink]{cleveref}
\usepackage{verbatim}
\usepackage{enumerate}
\usepackage{slantsc}

\newcommand{\R}{\mathbb{R}}

\global\long\def\R{\mathbb{R}}

\global\long\def\S{\mathbb{S}}

\global\long\def\L{\mathcal{L}}

\global\long\def\P{\mathcal{P}}

\global\long\def\F{\mathcal{F}}

\global\long\def\S{\mathcal{S}}

\newtheorem{theorem}{Theorem}[section]
\newtheorem*{namedtheorem}{\theoremname}
\newcommand{\theoremname}{testing}

\newtheorem{thm}[theorem]{Theorem}

\newtheorem{corollary}[theorem]{Corollary}

\newtheorem*{question*}{Question}

\theoremstyle{definition}
\newtheorem{definition}[theorem]{Definition}

\newtheorem{remark}[theorem]{Remark}

\theoremstyle{plain}

\title{A note on the largest bipartite subgraph in point-hyperplane incidence graphs}
\author{Thao~Do \thanks{Massachusetts Institute of Technology, Department of Mathematics. Email: thaodo@mit.edu} }

\begin{document}
\maketitle
\begin{abstract}
Given $m$ points and $n$ hyperplanes in $\R^d$, if there are many incidences, we expect to find a big cluster $K_{r,s}$ in their incidence graph. In \cite{apfelbaum2007large}, Apfelbaum and Sharir found lower and upper bounds for the largest size of $rs$, which only match in three dimensions. In this paper we close the gap in four and five dimensions, up to some logarithmic factors. 
\end{abstract}

\section{Introduction}
 Given a set $P$ of $m$ points and a set $Q$ of $n$ hyperplanes in $\R^d$, their \emph{incidence graph} $G(P,Q)$ is a bipartite graph with vertex set $P\cup Q$ and $(p,q)\in P\times Q$ forms an edge iff $p\in q$. It is proved in \cite{apfelbaum2007large} that if this graph does not contain $K_{r,s}$ as a subgraph for some fixed $r,s$, then it can have at most $O_d((mn)^{d/(d+1)}+m+n)$ edges. Here the notations $f=O_d(g)$  means there exists some constant $C$ that depends on $d$ such that $f\leq Cg$. 
 
 Conversely, when the graph has many edges, we expect to find a big $K_{r,s}$ subgraph. How big can $rs$ be in term of $m,n$ and the number of edges? To make it precise, we use the following definition.
\begin{definition}\label{def:rs}
Given a set $P$ of points and $Q$ of hyperplanes in $\R^d$, let $rs(P,Q)$ be the maximum size of a complete bipartite subgraph of its incidence graph, and $rs_d(m,n,I)$ be the minimum of this quantity over all choices of $m$ points and $n$ hyperplanes in $\R^d$ with $I$ incidences. To be precise: $$rs(P,Q):=\max\{rs: K_{r,s}\subset G(P,Q)\}$$
$$rs_d(m,n,I):=\min_{|P|=m,|Q|=n,\newline |G(P,Q)|=I} rs(P,Q).$$
\end{definition}

Apfelbaum and Sharir in \cite{apfelbaum2007large} proved that if $I=\Omega_d( mn^{1-\frac{1}{d-1}}+ nm^{1-\frac{1}{d-1}})$ then 
\begin{equation}\label{eq:ap-sh lower}
    rs_d(m,n,I)=\Omega_d\left( \left(\frac{I}{mn}\right)^{d-1} mn\right).
\end{equation}
Moreover, they show the following upper bound: if $I=\Omega_d((mn)^{1-\frac{1}{d-1}})$ then
\begin{equation}\label{eq:ap-sh upper}
    rs_d(m,n,I)=O_d\left( \left(\frac{I}{mn}\right)^{\frac{d+1}{2}} mn\right).
\end{equation}

These lower and upper bounds only match when $d=3$. In this paper we close the gap in four and five dimensions. In particular, we improve the lower bound to match with the upper bound up to some logarithmic factors.
 \begin{thm}\label{thm:4d}
When $d=4$, there exist constants $C_4$ and $C'_4$ such that if $I\geq C_4  (mn^{2/3}+nm^{3/5})$ then
$$rs_4(m,n,I)\geq C'_4 \left(\frac{I}{mn}\right)^{5/2} mn(\log mn)^{-4}.$$
\end{thm}
 \begin{thm}\label{thm:5d} When $d=5$, there exist constants $C_5$ and $C'_5$ such that if $I\geq C_5   (mn^{3/4}+nm^{2/3})$ then
$$rs_5(m,n,I)\geq C'_5 \left(\frac{I}{mn}\right)^{3} mn(\log mn)^{-10}.$$
\end{thm}

The main tool used to prove \cref{thm:4d} and \cref{thm:5d} is an incidence bound between points and \emph{nondegenerate} hyperplanes, which is reviewed in the next section. We then present the proof of \cref{thm:4d} and sketch the proof of \cref{thm:5d} in the subsequent sections. At the end we explain why our method does not work in six dimensions.
\subsection*{Acknowledgement} The author would like to thank Larry Guth  for many helpful conversations.

\section{Incidence with nondegenerate hyperplanes}
We use the following notations. Let $A$ and $B$ be two sets of geometric objects in $\R^d$. Their \emph{incidence graph} $G(A,B)$ is a bipartite graph on $A\times B$ where $(a,b)$ forms an edge iff $a\subset b$. The number of incidences between $A$ and $B$, denoted by $I(A,B)$, is the number of edges of this graph. In this paper, $A$ are either a set of points or a set of lines, and $B$ is a set of higher dimensional flats. Moreover, we sometimes use the notations $ f \gtrsim g$ instead of $f=\Omega(g)$ and $f\lesssim g$ instead of $f=O(g)$.

Given a set $\S$ of $m$ points in $\R^d$ and some $\beta\in(0,1)$, a hyperplane $H$ is \emph{$\beta$-nondegenerate} with respect to (w.r.t.) $\S$ if there does not exist a proper subflat $F\subset H$ that contains more than $\beta$ fraction of the number of points of $\S$ in $H$, i.e. $|F\cap \S|>\beta|H\cap \S|$. Otherwise, $H$ is \emph{$\beta$-degenerate}. Elekes and T\'oth proved the following incidence bound.

\begin{thm}[Elekes-T\'oth \cite{elekes2005incidences}]\label{thm:E-T}
If $\S$ is a set of $m$ points and $\mathcal{H}$ is a set of $n$ $\beta$-nondegenerate hyperplanes w.r.t. $\S$ (for any $0<\beta<1$)\footnote{Elekes-T\'oth actually proved this only for $\beta<\beta_d$ for some small $\beta_d$. It is later shown in \cite{do2016extending} that we can take $\beta_d=\frac{1}{d-1}$ and in \cite{lund2017two} that we can take $\beta_d=1$.} in $\R^{d}$ then there exists a constant $C_{\beta, d}$ such that 
\begin{equation}\label{eq:E-T}
I(\S,\mathcal{H})\leq C_{\beta,d}\left((mn)^{\frac{d}{d+1}}+mn^{1-\frac{1}{d-1}}\right).
\end{equation}
\end{thm}
This implies the maximum number of $\beta$-nondegenerate, $k$-rich (i.e. containing at least $k$ points of $\S$) hyperplanes is  
$O_{\beta,d} \left(\frac{m^{d+1}}{k^{d+2}}+\frac{m^{d-1}}{k^{d-1}}\right).$
Actually this is what Elekes-T\'oth proved. It is shown to be equivalent to \eqref{eq:E-T} in \cite{apfelbaum2007large}. 

Since lines and hyperplanes are dual to each other, we also have a dual version of the above result. Given a set $\mathcal{H}$ of n hyperplanes in $\R^d$, a point $p$ is $\beta$-nondengenerate with respect to $\mathcal{H}$ if there does not exist a line $\ell$ such that $\#\{H\in \mathcal{H}: \ell\subset H\}\geq \beta \#\{H\in \mathcal{H}: p\in H\}$.
\begin{corollary}\label{cor:dual E-T}
If $\mathcal{H}$ is a set of $n$ hyperplanes in $\R^d$ and $P$ is a set of $m$ $\beta$-nondegenerate points w.r.t. $\mathcal{H}$ then there exists a constant $C'_{\beta,d}$ such that
\begin{equation}\label{eq:E-T-dual}
I(\S,\mathcal{H})\leq C'_{\beta,d} \left((mn)^{\frac{d}{d+1}}+nm^{1-\frac{1}{d-1}}\right).
\end{equation}
Equivalently, given $n$ hyperplanes in $\R^d$, the number of $k$-rich, $\beta$-nondegenerate points is $O_{\beta,d} \left(\frac{n^{d+1}}{k^{d+2}}+\frac{n^{d-1}}{k^{d-1}}\right).$
\end{corollary}

\section{Proof in four dimensions}

We first outline our strategy. Let $\S$ be a set of $m$ points, $\mathcal{H}$ be a set of $n$ hyperplanes in $\R^4$. There are two ways to form a big $K_{r,s}$ in the incidence graph $G(\mathcal{H},\S)$: either a plane contains many points of $\S$ and belongs to many hyperplanes of $\mathcal{H}$, or a line does. By an averaging argument, we can assume each hyperplane is $\frac{I}{m}$-rich (i.e. contains at least $\frac{I}{m}$ points of $\S$). By \cref{thm:E-T}, the contribution from $\beta$-nondegenerate hyperplanes is negligible, so we can assume each hyperplane is $\beta$-degenerate, i.e. it contains some plane with at least $\beta$ of the total number of points, hence the plane is $\beta\frac{I}{m}$-rich. In this case, we say each hyperplane \emph{degenerates} to a rich plane. Either one of those planes belongs to many hyperplanes, which would form a big $K_{r,s}$, or we can find a subset $\P_i$ of planes such that $I(\S,\P_i)$ is large.  We repeat our argument: using the averaging argument and \cref{cor:dual E-T}, we can assume each point in $\S$ belongs to many planes in $\P_i$ and degenerates to a line. Either one of those lines contains many points, which then form a big $K_{r,s}$, or we can find a subset $\L_j$ of lines such that $I(\L_j,\P_i)$ is large. But after some transformation, this number is the same with the number of incidences between points and lines in $\R^2$ and hence cannot be too large by \cref{thm:E-T} for $d=2$, or equivalently, Szemer\'edi-Trotter's theorem in \cite{szemeredi1983extremal}. 

We now give the detailed proof.

\begin{proof}[Proof of \cref{thm:4d}] Assume $I\geq C_4(mn^{2/3}+nm^{3/5})$ for some big constant $C_4$ chosen later, but the incidence graph $G(\S, \mathcal{H})$ with $I$ edges contains no $K_{r,s}$ of size $rs\gtrsim \left(\frac{I}{mn}\right)^{5/2} mn(\log mn)^{-4}$. We follow several steps to derive a contradiction.
\\
\\
\noindent \textbf{Step 1:} We can assume each hyperplane is $\frac{I}{4n}$-rich and $\beta$-degenerate with respect to $\S$ for some $\beta>0$. 

Indeed,
remove all the hyperplanes that contain fewer than $\frac{I}{4n}$ points and the hyperplanes that is $\beta$-nondegenerate. The number of incidences from the non-rich hyperplanes is at most $n\frac{I}{4n}=\frac{I}{4}$. By \cref{thm:E-T}, the number of incidences from the $\beta$-nondegenerate hyperplanes  is at most $C_{\beta,4}((mn)^{4/5}+mn^{2/3})<\frac{C_4}{4}(mn^{2/3}+nm^{3/5})$ for big enough $C_4$. Indeed, this only fails if $(mn)^{4/5}\gtrsim  mn^{2/3}$ and $(mn)^{4/5}\gtrsim nm^{3/5}$, which is equivalent to $m\lesssim n^{2/3}$ and $n\lesssim m$, but they cannot happen at the same time for appropriate choices of constants. Therefore, after the removal, there remains at least $\frac{I}{2}$ incidences left. Assume there are $n_1$ hyperplanes left, where $n_1\leq n$. In fact, throughout the proof, we always use $n$ to upper bound $n_1$, so we can simply assume $n_1=n$.
\\
\\
\noindent \textbf{Step 2:} 
For each $\frac{I}{4n}$-rich $\beta$-degerenate hyperplane $H$, we can find a plane $P\subset H$ so that $|P|\geq \beta|H|\geq \frac{\beta I}{4n}$. Let $\mathcal{P}$ denote the set of these planes. We claim that no plane in $\mathcal{P}$ belongs to more than $s_0$ hyperplanes in $\mathcal{H}$ where 
\begin{equation}\label{eq:s}
s_0:=\frac{c_1 I^{3/2}}{m^{3/2}n^{1/2}(\log mn)^4}\end{equation} 

Indeed, assume there are $\frac{c_1I^{3/2}}{m^{3/2}n^{1/2}(\log mn)^4}$ hyperplanes that degenerate to a same plane for some constant $c_1$, then we have a configuration of $K_{r,s}$ where $$rs\geq \frac{\beta I}{4n}\cdot\frac{c_1I^{3/2}}{m^{3/2}n^{1/2}(\log mn)^4}\geq C'\left(\frac{I}{mn}\right)^{5/2} mn(\log mn)^{-4}$$
if we choose $C'<\frac{\beta c_1}{4}$. Contradiction.
\\
\\
\noindent \textbf{Step 3:} We use a dyadic decomposition to find a subset of planes with lots of incidences with $\S$.
Let $\mathcal{P}_j$  denote the set of all planes that is assigned to at least $2^j$ and at most $2^{j+1}$ hyperplanes where $j<\log s_0<\log n$ (here the logarithm is in base 2). We claim that there exists some $i$ such that
\begin{equation}\label{eq:I'}
    I':=I(\S,\P_i)> 4C_{\beta,3} (|\P_i||\S|)^{3/4}+|\P_i||\S|^{1/2}
\end{equation} where $\beta$ is the same as before, and $C_{\beta,3}$ is defined in \cref{thm:E-T}.  

Indeed, first notice that the contribution to incidences from the planes must be at least $\beta$-fraction the number of incidences from the $\beta$-degenerate hyperplanes, which implies
$\sum_{j=0}^{\log s_0} 2^{j+1} I(\S,\P_j)\geq \frac{\beta}{4}I.$
Hence there must exist some $i$ such that
$$
2^{i_0+1}I(\S,\P_{i})\geq \frac{\beta I}{4\log s_0}.
$$
On the other hand, since each hyperplane is assigned to exactly one plane, we have
 $\sum_{j=0}^{\log s_0} 2^{j} |\P_j|\leq n_2\leq n.$
 As a consequence, $|P_i|\leq \frac{n}{2^i}$.
From \eqref{eq:s}, we have $2^i\leq s_0\leq \frac{c_1 I^{3/2}}{m^{3/2}n^{1/2}(\log mn)^4}$. Assume \eqref{eq:I'} fails, then
\begin{align*}
\frac{\beta I}{4\log s_0} & \leq 2^{i+1} I'\\
&\leq 2^{i+1} 4C_{\beta,3} \left(|\P_i||\S|)^{3/4}+|\P_i||\S|^{1/2}\right) \\
& \leq 8C_{\beta,3} 2^{i}\left(\left(\frac{nm}{2^i}\right)^{3/4}+\frac{n}{2^i}m^{1/2}\right) \\ & \leq 8C_{\beta,3}\left((mn)^{3/4} \left(\frac{c_1 I^{3/2}}{m^{3/2}n^{1/2}(\log mn)^4}\right)^{1/4}+ nm^{1/2}\right)
\end{align*}
Since $I\geq Cnm^{3/5}\gg nm^{1/2}\log s_0$, the first term in the right hand side must be at least $\frac{\beta I}{8\log s_0}$. Rearranging we get 
$$I^{5/8}\leq c_3 m^{3/8}n^{5/8}\frac{\log s_0}{\log mn}$$
where $c_3$ depends on $\beta, C_{\beta,3}$ and $c_1$. However, we can choose $\beta$ and $c_1$ small enough and $C_4$ big enough so that $c_3^{8/5}<C_4$ and hence this contradicts with $I\geq C_4 nm^{3/5}$. So \eqref{eq:I'}
must hold.
\\
\\
\noindent\textbf{Step 4:} Since the bound in \eqref{eq:I'} is the same with that in \cref{cor:dual E-T}, we can use a similar argument with Step 1 to assume each point in $S$ is $\frac{I'}{4m}$-rich (i.e. belongs to at least $\frac{I'}{4m}$ planes in $\P_i$), and is $\beta$-degenerate w.r.t. $\P_i$ (in the sense defined before \cref{cor:dual E-T}. Each such point degenerates to a line that is $\beta\frac{I'}{4m}$-rich. Let $\L$ denote the set of all these lines. We claim that no line in $\L$ contains more than $r_0$ points where 

\begin{equation}\label{eq:r}
r_0:=\frac{c_5 I^{3/2}}{m^{1/2}n^{3/2}(\log mn)^3}
\end{equation} 


Indeed, each line in $\L$ belongs to at least $\frac{\beta I'}{4m}$ planes in $\P_i$,  and thus belongs to at least $\frac{\beta I'2^i}{4m}\geq \frac{\beta I}{4m\log s_0}$  hyperplanes in $\mathcal{H}$ because each plane in $\P_i$ belongs to at least $2^i$ hyperplanes. 
If there are $\frac{c_5 I^{3/2}}{m^{1/2}n^{3/2}(\log mn)^3}$ points that degenerates (or belongs) to a same line for some constant $c_5$, then we have a configuration of $K_{r,s}$ where 
$$rs\geq \frac{\beta I}{4m\log s_0}\cdot\frac{c_5 I^{3/2}}{m^{1/2}n^{3/2}(\log mn)^3}\geq C'_4\left(\frac{I}{mn}\right)^{5/2} mn(\log mn)^{-4}$$
for small enough $C'_4$. Contradiction.
\\
\\
\noindent\textbf{Step 5:} Similar to Step 3, we use a dyadic decomposition to find a subset of lines in $\L$ that has many incidences with $\P$. Here we say a line $\ell$ is incident to a plane $P$ if $\ell\subset P$. 
Let $\mathcal{L}_k$  denote the set of all lines that contain at least $2^k$ and at most $2^{k+1}$ points where $k<\log r_0<\log m$. Note that here for a line we only consider the points that degenerate to that line. We claim there must exist some $j$ such that 
\begin{equation}\label{eq:I''}
    I'':=I(\L_j,\P_i)\geq C_{\beta, 2}\left( |\P_i|^{2/3}|\L_j|^{2/3}+|\P_i|+|\L_j|\right)
\end{equation}

Indeed, first notice that the contribution to incidences from the lines must be at least $\beta$-fraction, which implies
$\sum_{k=0}^{\log r_0} 2^{k+1} I(\L_k,\P_i)\geq \frac{\beta}{4}I'.$
Hence there must exist some $j$ such that 
$$2^{j+1}I(\L_j,\P_i)\geq \frac{\beta I'}{4\log r_0 }\geq \frac{\beta^2 I}{16\log s_0\log r_0 2^{i}}.
$$
On the other hand, since each point is assigned to exactly one line  we have
$ \sum_{k=0}^{\log r_0} 2^{k} |\L_k|\leq m.$ As a consequence, $|\L_j|\leq\frac{m}{2^j}$. From \eqref{eq:r}, $2^j\leq r_0\leq \frac{c_5 I^{3/2}}{m^{1/2}n^{3/2}(\log mn)^3}$. Recall $|\P_i|\leq\frac{n}{2^i}$ and $2^i\leq s_0$. Assume \eqref{eq:I''} fails, then there exists some constant $c_6$ such that

\begin{align*}
I & \leq 32\beta^{-1}(\log r_0\log s_0) 2^{i+j}I(P_i, S_j)\\
& \leq c_6 (\log m\log n)\left( (s_0 r_0)^{1/3}(mn)^{2/3}+nr_0+ ms_0\right)\\
&\leq c_6(\log m\log n)\left( (mn)^{2/3}\left(\frac{c_1 I^{3/2}}{m^{3/2}n^{1/2}(\log mn)^4}\frac{c_5 I^{3/2}}{m^{1/2}n^{3/2}(\log mn)^3}\right)^{1/3}+m \frac{c_1 I^{3/2}}{m^{3/2}n^{1/2}(\log mn)^4}+n\frac{c_5 I^{3/2}}{m^{1/2}n^{3/2}(\log mn)^3}\right)\\
&\leq I\left[\frac{c_1c_5c_6\log m\log n}{(\log mn)^7}+\left(\frac{I}{mn}\right)^{1/2}\left(\frac{c_1\log m\log n}{(\log mn)^4}+\frac{c_5\log m\log n}{(\log mn)^3}\right)\right]\\
&\ll I.
\end{align*}
This contradiction implies \eqref{eq:I''} must hold. 
\\
\\
\textbf{Step 6:} We claim that \eqref{eq:I''} violates \cref{thm:E-T} in two dimensions, or Szemer\'edi-Trotter's theorem. 

Indeed, project the set of planes $\P_i$ and the set of lines $\L_j$ to a generic three dimensional subspace, then intersect them with  generic plane $\Pi$ within this subspace. After this transformation, $\P_i$ becomes a set of lines $P^*$ and  $\L_j$ becomes a set of points $L^*$ in $\Pi$. We have $I(P^*, L^*)=I(\L_j,\P_i)\gtrsim |P^*|^{2/3}|L^*|^{2/3}+|P^*|+|L^*|$. This violation finishes our proof.
 \end{proof}

\section{Sketch of proof in five dimensions}
The proof method is the same with that in four dimensions, but the exponents are different and the method is repeated one more time.
In particular, we unwrap in three layers: hyperplanes degenerate to 3-flats, points degenerate to lines, and 3-flats degenerate to planes. At each layer, either we can find a big $K_{r,s}$, or the number of incidences remain larger than the nondegenerate bound in \cref{thm:E-T}, and we can keep unwrapping. The detailed proof  is quite similar to that in the four dimensions case, so we only give an outline here. For simplicity, we ignore all the constants and logarithmic factors.

\begin{proof}[Proof's sketch of \cref{thm:5d}] Prove by contradiction. Let $\S$ denote the set of $m$ points and $\mathcal{H}$ denote the set of $n$ hyperplanes in $\R^5$.
Assume $I(\S,\mathcal{H})\geq C_5(mn^{3/4}+nm^{2/3})$ but their incidence graph does not contain any $K_{rs}$ where $rs\gtrsim \left(\frac{I}{mn}\right)^3 mn(\log mn)^{-10}$. 

\begin{itemize}
\item[Step 1] We can assume every hyperplane is $\frac{I}{n}$-rich, and $\beta$-degenerate with respect to $\S$ for some $\beta>0$.
\item[Step 2] For each such hyperplane $H$, we can find a 3-dimensional flat (or a 3-flat) $F$ such that $F\subset H$ and $|F\cap \S|\geq \beta|H\cap \S|\geq\frac{\beta I}{n}$. Let $\F$ denote the set of these 3-flats. We show that no flat in $\F$
belong to more than $s_0$ hyperplanes where $s_0\lesssim \frac{I^2}{m^2 n}$.

\item[Step 3] Let $\F_j$ denote the set of all 3-flats in $\F$ that is assigned to at least $2^j$ and at most $2^{j+1}$ hyperplanes where $j\leq \log s_0<\log n$. We show that there exists an $i$ such that 
$$I':=I(\F_i, \S)\gtrsim (|\F_i||\S|)^{4/5}+|\F_i||\S|^{2/3}.$$
Indeed, assume otherwise. Using $I'\gtrsim 2^i I$, $|\F_i|\leq \frac{n}{2^i}$ and $2^i\leq s_0\lesssim \frac{I^2}{m^2n}$, we have
$$I\lesssim 2^i I'\lesssim  2^i\left[\left(\frac{nm}{2^i}\right)^{4/5}+\frac{n}{2^i}m^{2/3}\right]
\lesssim (mn)^{4/5}\left(\frac{I^2}{m^2n}\right)^{1/5}+nm^{2/3}$$
which cannot happen given our condition $I\gtrsim mn^{3/4}+nm^{2/3}$.
\item[Step 4] Since $I'$ is large, using \cref{cor:dual E-T}, we can assume each point in $\S$ is $\frac{I'}{m}$-rich (i.e. belongs to at least $\frac{I'}{m}$ flats in $\F_i$, and is $\beta$-degenerate w.r.t. $\F_i$.
Each such point degenerates to a $\frac{\beta I'}{m}$-rich line. Let $\L$ denote that set of these lines. Then no line in $\L$ can contain more than $r_0$ points where $r_0\lesssim \frac{I^2}{mn^2}$.
\item [Step 5] We use a dyadic decomposition to find a subset of lines with many incidences with $\F_i$. Let $\L_k$ denote the set of all lines in $\L$ that contain more than $2^k$ and fewer than $2^{k+1}$ points. Then 
there exists a $j$ such that 
$$I'':=I(\F_i, \L_j)\gtrsim |\F_i|^{3/4}|\L_j|^{3/4}+|\F_i||\L_j|^{1/2}.$$
Indeed, assume otherwise. Using $I''\gtrsim I' /2^j\gtrsim I /2^{i+j}$ $|\F_i|\leq\frac{n}{2^i}$, $|\L_j|\leq\frac{m}{2^j}$, $2^i\leq s_0\lesssim \frac{I^2}{m^2n}$ and $2^j\leq r_0\lesssim\frac{I^2}{mn^2}$, we have
$$I\lesssim 2^{i+j} I''\lesssim 2^{i+j}\left[\left(\frac{mn}{2^{i+j}}\right)^{3/4}+\frac{n}{2^i}\left(\frac{m}{2^j}\right)^{1/2}\right]\lesssim (mn)^{3/4} \left(\frac{I^2}{m^2n}\frac{I^2}{mn^2}\right)^{1/4}+nm^{1/2}\left(\frac{I^2}{mn^2}\right)^{1/2}=2I$$
This cannot happen with an appropriate choice of constants and logarithmic factors.
\item [Step 6] Turn $I(\F_i, \L_j)$ into point-plane incidences in $\R^3$ by some transformation.  This means we can assume each 3-flats in $\F_i$ degenerate to a plane. Let $\P$ denote the set of all such planes. Then no plane belongs to more than $t_0$ flats in $\F_i$ where $t_0\lesssim \frac{I^2}{m^2n}$
\item[Step 7] Using a dyadic decomposition, there exists some subset $\P_k$ of planes, each belongs to at least $2^k$ and at most $2^{k+1}$ planes in $\F_i$ such that
$I''':=I(\L_j,\P_k)\gtrsim |\P_k|^{2/3}|L_j|^{2/3}+|\P_k|+|\L_j|.$
\item[Step 8] Turn $I'''$ into point-line incidences, which leads to a violation with Szemer\'edi-Trotter's theorem. This finishes our proof.
\end{itemize}
\end{proof}
\begin{remark}
Our argument does not work in six dimensions and higher because when we write down the details of step 5 in the above outline, the exponents no longer match and we do not get a contradiction. 
\end{remark}

\bibliographystyle{abbrv}
\bibliography{large_subgraph.bib}

\begin{thebibliography}{1}

\bibitem{apfelbaum2007large}
R.~Apfelbaum and M.~Sharir.
\newblock Large complete bipartite subgraphs in incidence graphs of points and
  hyperplanes.
\newblock {\em SIAM Journal on Discrete Mathematics}, 21(3):707--725, 2007.

\bibitem{do2016extending}
T.~Do.
\newblock Extending {E}rd{\H o}s -{B}eck's theorem to higher dimensions.
\newblock {\em arXiv:1607.00048}, 2016.

\bibitem{elekes2005incidences}
G.~Elekes and C.~D. T{\'o}th.
\newblock Incidences of not-too-degenerate hyperplanes.
\newblock In {\em Proceedings of the twenty-first annual symposium on
  Computational geometry}, pages 16--21. ACM, 2005.

\bibitem{lund2017two}
B.~Lund.
\newblock Two theorems on point-flat incidences.
\newblock {\em arXiv:1708.00039}, 2017.

\bibitem{szemeredi1983extremal}
E.~Szemer{\'e}di and W.~T. Trotter.
\newblock Extremal problems in discrete geometry.
\newblock {\em Combinatorica}, 3(3-4):381--392, 1983.

\end{thebibliography}
\end{document}